\documentclass[11pt]{article}%
\usepackage{ifpdf}
\usepackage{amsmath}
\usepackage{amsfonts}
\usepackage{amssymb}
\usepackage{verbatim}
\usepackage{psfrag}
\usepackage{graphicx}
\usepackage[T1]{fontenc}
\usepackage{lmodern}%
\ifpdf
\usepackage{garamond}
\fi
\providecommand{\U}[1]{\protect\rule{.1in}{.1in}}
\providecommand{\U}[1]{\protect\rule{.1in}{.1in}}
\newtheorem{theorem}{Theorem}

\newtheorem{corollary}[theorem]{Corollary}

\newtheorem{definition}[theorem]{Definition}

\newtheorem{lemma}[theorem]{Lemma}

\newtheorem{proposition}[theorem]{Proposition}

\newenvironment{proof}[1][Proof]{\noindent\textbf{#1.} }{\ \rule{0.5em}{0.5em}}

\begin{document}

\ifpdf
\garamond
\fi
\title{Semigroups over Local fields}
\author{Marcelo Firer (Universidade Estadual de Campinas) \\ and Daniel Miranda (Universidade Federal do ABC)}
\date{}
\maketitle
\begin{abstract}
Let $G$ be a $1$-connected, almost-simple Lie group over a local field and
$\mathcal{S}$ a subsemigroup of $G$ with non-empty interior. The action of the
regular hyperbolic elements in the interior of $\mathcal{S}$ on the flag
manifold $G/P$ and on the associated Euclidean building allows us to prove
that the invariant control set exists and is unique. We also provide a
characterization of the set of transitivity of the control sets: its elements
are the fixed points of type w for a regular hyperbolic isometry, where w is
an element of the Weyl group of $G$. Thus, for each w in W there is a control
set $D_{w}$ and $W(\mathcal{S})$ the subgroup of the Weyl group such that the
control set $D_{w}$ coincides with the invariant control set $D_{1}$ is a Weyl
subgroup of $W$. We conclude by showing that the control sets are
parameterized by the lateral classes $W(S)\backslash W.$

\end{abstract}
\tableofcontents


\newpage

\newpage

\section{Semigroups and Control Sets}

Our purpose in this text is to study a particular class of semigroups: namely,
the semigroup of a Lie group over a local field. This task will be tackled in
an essentially geometric way, by means of a detailed study of the action of
the semigroup in a suitable space. This approach leads to a rich theory that
relates the semigroup to the control sets. We begin with some basic definitions:

Let $X$ be a topological space and $C(X)$ the group of homeomorphisms of
$X$\textbf{,} a semigroup (of homeomorphism) is a set $\mathcal{S}\subset
C(X)$ closed under composition.

Let $\mathcal{S}$ be a semigroup of homeomorphism of $X$ and $x$ a point in
$X.$ The set $\mathcal{S}x=\{gx\in X:g\in\mathcal{S\}}$ is called the
\textbf{orbit} of $x$ by the action of $\mathcal{S}.$ We refer to the closure
$\operatorname*{cl}(\mathcal{S}x)$ of $Sx$ in $X$ as the \textbf{approximate
orbit}. \ A subset of $X$ is $\mathcal{S}$-invariant if $\mathcal{S}D\subset
D.$

We say that the action of $\mathcal{S}$ on a set $A\subseteq X$ is
\textbf{approximately transitive} if $A\subseteq\operatorname*{cl}%
(\mathcal{S}x)$ for all $x\in A$.

\begin{definition}
A \textbf{control set} for $S$ in $X$ is a subset $D\subset X$ which fulfills
the following conditions :

\begin{enumerate}
\item $D\subset\operatorname*{cl}(\mathcal{S}x)$ for all $x\in D$;

\item $D$ is maximal with respect to the first property;

\item $\operatorname*{int}D\neq\emptyset$
\end{enumerate}
\end{definition}

The first two items of the previous definition of control set say that these
sets are maximal among all sets where the action of $\mathcal{S}$ is
approximately transitive.

We introduce the following partial ordering between control sets: $D_{1}$ is
smaller than $D_{2}$ if there exist $x\in D_{1}$ and $s\in\mathcal{S}$ such
that $sx\in D_{2}$. A maximal control set with respect to this ordering is
called an \textbf{invariant control set} ,which is clearly $\mathcal{S}$-invariant.

We define the \textbf{set of transitivity} as%
\[
D_{0}=\{x\in D:x\in\operatorname*{int}(\mathcal{S}^{-1}x)\}
\]

\begin{proposition}
Let $\mathcal{S}$ be a semigroup of homeomorphisms with non-empty interior,
$D$ a control set for $\mathcal{S}$ in $X$ and $D_{0}$ the set of transitivity
of $D$. Then the following hold:

\begin{enumerate}
\item $D_{0}=\operatorname*{int}(\mathcal{S}x)$ for all $x\in D_{0}$

\item $cl(D_{0})=D$

\item For any $x,y\in D_{0}$\textbf{\ }there is\textbf{\ }$s\in S$ with $sx=y
$
\end{enumerate}
\end{proposition}

We point out that item (3) above is the reason why $D_{0}$ is called the set
of transitivity of $D$.

When the topological space $X$ is compact, an application of Zorn's Lemma
yields the existence of an invariant control set $D$.

In the very special case in which $S\subset G$\textbf{,} with $G$ a
semi-simple Lie group, San Martin has developed a powerful theory relating the
semigroups and the control sets in the flag manifold $G/P$. Among others
results San Martin proved that if $\operatorname*{int}(\mathcal{S)\neq
\emptyset}$ then there is a unique control set in $G/P$\textbf{,} and that for
each element of the Weyl group $w\in W$ there exist a control set $D_{w}$ on
maximal flag $G/P$, whose elements of the set of transitivity are fixed points
of type $w$\textbf{,} for some regular $h$\textbf{.} For more information and
another results see \cite{San1}, \cite{San2}, \cite{San3} and \cite{ST}.

We extend these results for algebraic groups over a local field.

\vspace{.5cm}

In sections \ref{2}, \ref{3} and \ref{4}, we introduce the basic concepts and
definitions used in this work. Since those are quite extense, the
presentations is somehow schematic form. In each of these sections we point
our basic bibliographic source but warn the reader that, due to notations used
in different contexts, the notations used in this work may differ slightly
from the usually adopted in the indicated bibliographical resources.

\section{Buildings\label{2}}

The basic bibliographical resource we used for this section are \cite{Bro} and
\cite{Ga}.

Let $I$ be a finite set of indexes. A \textbf{Coxeter matrix} $M=(m_{ij})$ is
a square matrix with values in $\mathbb{N}\cup\{\infty\}$ such that $m_{ij}=1$
if and only if $i=j$ and $m_{ij}=m_{ji}.$

A Coxeter matrix defines a \textbf{Coxeter group} of type $M$, the group
$W(M)$ defined by the presentation
\[
W(M)=\left\langle r_{i}:r_{i}^{2}=(r_{i}r_{j})^{m_{ij}}=1,\forall i,j\in
I\right\rangle
\]

To simplify the notation we denote $W(M)$ by $W$ wherever it is clear which
Coxeter Matrix is been referred to.

We can see identify the set of indexes $I$ with the the set $\{r_{i}:i\in I\}
$ of order two generators of $W$ and call the pair $(W,I)$ a \textbf{Coxeter
system}.

A\textbf{\ special subgroup} of a Coxeter group $W$ with generators $J$ is a
subgroup $W(J)$ generated by a subset $J\subseteq I$. A \textbf{special coset}
is a coset $wW(J)$, determined by a special subgroup $W(J)$.

We describe now the chamber complex system, the \textbf{Coxeter complex}
associated to a Coxeter system, which plays a crucial r\^{o}le\texttt{\ }in
this work. Let $\mathcal{V}$ be a poset with the usual inclusion order. Let
$\Delta$ be a family of finite subsets of $\mathcal{V}$ containing every
singleton $\{v\}\subseteq V$ and satisfying the condition that, if $A\in
\Delta$ and $B\subseteq A$ then $B\in\Delta$. Such a pair $(\mathcal{V}%
,\Delta)$ is called a \textbf{simplicial complex} and $\Delta$ a
\textbf{family of simplices}. Given simplices $B\subseteq A\in\Delta$ we say
that $B$ is a \textbf{face of} $A$.

The cardinality $r$ of a simplex $A$ is said to be its \textbf{rank} and $r-1$
is called the \textbf{dimension} of $A$. Two simplices $A,B$ in a simplicial
complex $\Delta$ are \textbf{adjacent} if they have a codimension $1$ face. A
\textbf{gallery} is a sequence of maximal simplices in which any two
consecutive ones are adjacent.

We say that $\Delta$ is a \textbf{chamber complex} if all maximal simplices
have the same dimension and any two can be connected by a gallery. The maximal
simplices are called chambers. We say that a \textbf{chamber complex }$\Delta$
is \textbf{labelled } by a set $I,$ if there is a surjective function
$:V\rightarrow I$ \ that restricts to a bijection in each chamber.

Given a Coxeter group $W(I)$, we consider the set of all special cosets
$\Sigma=\Sigma(W,I)$. We define on $\Sigma$ a partial order inverse to the
inclusion: $W(J)\leq W(L)$ if $L\subseteq J$. The posets $(\Sigma,\leq)$ and
$(\mathcal{P}(I),\subseteq)$ are isomomorphic and thus a simplicial structure
can be induced on $\Sigma$, thus called the \textbf{Coxeter complex} of
$(W,I)$. The maximal elements in $\Sigma(W,I)$, are called chambers. They are
precisely the minimal special cosets, i.e., the sets with exactly one element
$w$ with $w\in W.$ The adjacency relation is determined as follows: $w_{1}%
\sim_{i}w_{2}$ iff $w_{1}=r_{i}w_{2}$ with $r_{i}\in I,$. We remark that every
Coxeter complex can be labelled by $I$, the set of generators of $W.$ When
considering the usual action of $W$ in the set of its cosets by left product,
we have that:

\begin{proposition}
\cite{Hu}

\begin{enumerate}
\item The group $W$ acts on $\Sigma(W,I)$ by type-preserving automorphisms.

\item The action of the group $W$ is transitive on the collection of simplices
of a given type.
\end{enumerate}
\end{proposition}

\begin{definition}
\label{definicaoedificio}A \textbf{building} is a simplicial complex $\Delta,$
together with a family $\mathcal{A}$ of subcomplexes (called the set of
apartments), satisfying the following axioms:
\end{definition}

\begin{enumerate}
\item Each apartment $\Sigma$ is a Coxeter complex.

\item For any two simplices $A,B\in\Delta$, there is an apartment $\Sigma$
containing both of them.

\item If $\Sigma$ and $\Sigma^{\prime}$ are two apartments containing $A$ and
$B$, there is an isomorphism $\Sigma$ $\rightarrow\Sigma^{\prime}$ fixing $A$
and $B$ pointwise.

\item Every codimension $1$ simplex is a face of at least three chambers
(thickness axiom).
\end{enumerate}

It is clear, by Axiom $B2$, that every two apartments are isomorphic. It also
implies that two maximal simplices $A$ and $B$ have the same dimension and can
be connected by a gallery in the apartment $\Sigma$ containing both of them,
and so $\Delta$ is a chamber complex. Furthermore $\Delta$ is a labelled
chamber complex labeled by $I$, the set of generators of the Coxeter group of
any given apartment $\Sigma$. The isomorphism of apartments $\Sigma
\rightarrow\Sigma^{\prime}$ postulated in Axiom $B2$ can be regarded as a
label-preserving isomorphism.

Any collection $\mathcal{A}$ of apartments $\Sigma$ satisfying the former
axioms is called a \textbf{system of apartments} for $\Delta$. It's well known
that the apartments are convex, i.e., given two chambers of $\Sigma$, then
every minimal gallery of $\Delta$ connecting these chambers is contained in
$\Sigma$.

\bigskip

Let $\Delta$ be a building, $A$ a system of apartments for $\Delta$ and $G$ a
group acting on $\Delta$ by simplicial label-preserving automorphisms that
leave $A$ invariant, i.e., if $\Sigma\in A$ then $g\Sigma\in A$.

We say that the action of $G$\textbf{\ }is\textbf{\ strongly transitive }if
$G$ acts transitively in the pairs $(\Sigma,C),$ with $C\in\Sigma$, in other
words, if $G$ is transitive in the set of apartments and the stabilizer of an
apartment is transitive on the chambers of $\Sigma$.

Henceforward we assume that $G$ acts strongly transitively in $\Delta,$ and
chose $\Sigma_{0}\in\mathcal{A}$ and $C_{0}\in\Sigma_{0}$, called the
\textbf{fundamental apartment} and the \textbf{fundamental chamber}, respectively.

The following subgroups of $G$ is of particular interest:%
\begin{align}
B  &  =\{g\in G:gC_{0}=C_{0}\}\label{bn}\\
N  &  =\{g\in G:g\Sigma_{0}=\Sigma_{0}\}\\
T  &  =\{g\in G:\ g\text{ fixes }\Sigma_{0}\text{ pointwise}\}\\
W  &  =N/T
\end{align}





Given $J\subseteq I$ we define $P_{J}\subseteq G$ as the subgroup of $G$
generated by $<B,W(J)>$, where $W(J)$ is the special subgroup with generators
indexed by $J$. A \textbf{face of type} $J$ is a simplex of $\Delta$ with
vertices labelled by $J$.

The following proposition, although not directly cited in the continuation is
crucial in many of the subsequents results in this section.

\begin{proposition}
Given $J^{\prime}\subset J$ and a face $A$ of type $J-J^{\prime}$, then the
stabilizer $P_{J^{\prime}}$ of $A$ satisfies:%
\[
P_{J^{\prime}}={\bigcup_{w\in W^{\prime}}}BwB.
\]

\end{proposition}

In particular, we have the \textbf{Bruhat Decomposition} of $G$:
\[
G={\coprod_{w\in W}}BwB
\]

As a consequence of the Bruhat Decomposition and the Axiom of Thickness, we
have that $BwB\circ BjB\subseteq BwB\cup BwjB$ and $jBj\varsubsetneq B$ for
every $j\in J$. All those properties are stated as postulates in the
definition of $BN$-pairs.

\begin{definition}
A\textbf{Tits system }or\textbf{\ $BN$-pair }is a group $G$ with two subgroups
$B,N$ satisfying:

\begin{description}
\item[BN0.] $B$ and $N$\ together generate $G$

\item[BN1.] $W=N/T$ is a Coxeter group with generators $J=\{j_{1},\ldots
j_{n}\}$, where $T:=B\cap N\vartriangleleft N$

\item[BN2.] $BjB\circ BwB\subset BwB\cup BjwB,$ for all $w\in W$ and $j\in J$

\item[BN3.] $jBj\neq B$ for all $j\in J$
\end{description}
\end{definition}

The relation between the structure of buildings and $BN$-pairs is the Tits' Theorem:

\begin{theorem}
Let $\Delta$\ be a building where a group $G$\ acts strongly transitively, and
let $B,C$\ and $N$\ be defined as in \ref{bn}. The pair $(B,N)$\ is then a
Tits system. Conversely, every Tits system $(B,N)$\ in a group $G$\ defines a
building in which the chambers are the cosets of $B$\ and the equivalence
relation is given by $gB\overset{j}{\sim}g^{\prime}B\Longleftrightarrow
gP_{j}=g^{\prime}P_{j}$. Finally, the action of $G$\ is strongly transitive
and $N$\ stabilizes an apartment.
\end{theorem}

\section{Geometric Realization of Buildings}

\label{3}

When the Coxeter group is an isometry group of a space we can give a geometric
interpretation to the Coxeter Complex and the associated building.

\subsection{Spherical Buildings}

Let $V$ be a finite-dimensional real vector space with an inner product,
$H\subset V$ an hyperplane and $\alpha$ an unit normal vector to $H$. The
reflection with respect to $H$ is the linear transformation defined as
$s_{H}(v)=v-\left\langle v,\alpha\right\rangle \alpha$\textbf{.}

A finite group $W$\ of linear transformations of $V$\ is called a finite
reflection group if it is generated by reflections $s_{H}$, where $H$ ranges
over a set $H$ of hyperplanes. In other words, $W$ is a discrete subgroup of
the \textbf{orthogonal group} $O(V)$\ which is generated by reflections.

We say that a finite reflection group $W$ is \textbf{essential }relative to
$V$ if $W$ acts on $V$ with no nonzero fixed points. Given a finite reflection
group $W$ and $V_{0}$ its space of fixed points, It is clear that any subgroup
$W$ stabilizes $V_{0}$ and its orthogonal complement $V_{0}^{\perp}$ and $W$
is essential relative to $V_{0}^{\perp}\,$.

A finite reflection group $\left(  W,V\right)  $ is called \textbf{reducible}
if $V$ decomposes as $V^{\prime}\oplus V^{\prime\prime}$ with $V^{\prime}$ and
$V"$ proper, $W$-invariant subspaces, i.e., $W(V^{\prime})\subset V^{\prime}$
and $W(V^{\prime\prime})\subset V^{\prime\prime}$\textbf{.}

Let $\mathcal{F}=\left\{  a_{H}=\left\langle \,\,,\alpha\right\rangle
:H\in\mathcal{H}\right\}  $ the set of all linear functionals associated with
$H$ $\in\mathcal{H}$. \ To each function $\phi:\mathcal{F}\rightarrow\left\{
+1,-1\right\}  $, we associate the set:%
\[
C_{\phi}=\{x\in V:\phi(x)\cdot f(x)>0,\forall f\in\mathcal{F\}}%
\]

The \textbf{chambers} are the non-empty $C_{\phi}$\textbf{.} The chambers form
a partition of $V\backslash\{H\in\mathcal{H}\}$ into disjoint convex cones.

A set $B$ is called a face of $A$ if its the intersection of $C_{\phi}$ with a
subspace of the form $S_{\mathcal{T}}=\{x\in X:f(x)=0,$ $\forall
f\in\mathcal{T}\subset\mathcal{H}\}$\textbf{.} In this case we write $B\leq
A$\textbf{.} A \textbf{wall} is a maximal proper face.

A Coxeter complex $\Sigma$ is called \textbf{spherical} if it is isomorphic to
the complex associated to a finite reflection group. Considering the
intersection of apartments, chambers, walls and faces with the unit sphere of
$V$ we get the geometric realization of the spherical Coxeter Complex. A
building $\Delta_{S}$ is called spherical if its apartments are spherical.

The diameter of a spherical building $\Delta_{S}$ is finite, and equal to the
diameter of any apartment. Two chambers $C,C^{\prime}$ in a spherical building
$\Delta_{S}$ are said to be opposite if $d(C,C^{\prime})=\operatorname*{diam}%
(A)$\textbf{.}

\subsection{Euclidean Buildings}

Let $V$ be a finite-dimensional real vector space with an inner product. The
isometry group of $V$ is the semidirect product $O(V)\ltimes V$\textbf{.} An
affine reflection group $W$ is a discrete subgroup of $O(V)\ltimes V$ that is
generated by reflections in affine hyperplanes.

We define $\mathcal{H}=\{H:s_{H}\in W\}$ and $\mathcal{F}=\{a_{H}%
:H\in\mathcal{H}\}$\textbf{.} We can then define the chambers for the affine
reflection group in a similar way as we did for the finite reflection groups.
\ In the affine case the chambers are bounded and they are the open Dirichlet
domains for the action of $W$ in $\mathbb{R}^{n}$.

Henceforth, we let $C_{0}$ be a fixed chamber (referred to as a
\textbf{fundamental chamber}) and let $J$ be the reflections in the walls of
$C_{0}$\textbf{.}

\begin{proposition}
\begin{enumerate}

\item The action of $W$ is simply transitive in the chambers.

\item $W$ is generated by $J$.

\item ($W,J)$ is a Coxeter Complex.

\item There exists $x\in V$ such that the stabilizer $W_{x}$ is isomorphic to
a finite Coxeter group $\overline{W}$.

\item $W\simeq\mathbb{Z}^{n}\ltimes$ $\overline{W}$\textbf{,} with $n$ being
the dimension of $V$\textbf{.}
\end{enumerate}
\end{proposition}

A point $x\in V$ such that $W_{x}$ is isomorphic to $\overline{W}$ is called
\textbf{special point}.

An abstract Coxeter complex $\Sigma$ is called Euclidean if it is isomorphic
to the complex $\left\vert \Sigma\right\vert $ associated to an affine
reflection group. The complex $\left\vert \Sigma\right\vert $ is called the
\textbf{geometric realization} of the complex $\Sigma$. We can choose a norm
in the vector space $V$ such that every chamber has diameter $1$. This norm is
called the \textbf{canonical norm} and with it every simplicial isomorphism
$\phi:\Sigma\rightarrow\Sigma^{\prime}$ induces an isometry $\left\vert
\phi\right\vert $ between the corresponding geometric realizations,
$\left\vert \phi\right\vert :\left\vert \Sigma\right\vert \rightarrow
\left\vert \Sigma\right\vert ^{\prime}$.

Let $\left\vert \Sigma\right\vert $ be the geometric realization of a
Euclidean Coxeter complex, and let $\mathcal{H}$ be the associated set of
hyperplanes in $V$. Fix $x\in V$ and let $\overline{\mathcal{H}}$ be the set
of hyperplanes through $x$ and parallel to some element of $\mathcal{H}$. The
set $\overline{H}$ is finite. Let $\mathcal{F}=\left\{  a_{H}=\left\langle
\,\,,\alpha\right\rangle :H\in\overline{\mathcal{H}}\right\}  $ the set of all
linear functional associated with $H$ $\in\overline{\mathcal{H}}$. \ Given a
function $\phi:\mathcal{F}\rightarrow\left\{  +1,-1\right\}  $, we associate
the set:%
\[
\mathcal{A}_{\phi}=\{y\in V:\phi(x)f(y-x)>0,\forall f\in\mathcal{F\}}%
\]

The \textbf{sectors} are the non-empty $\mathcal{A}_{\phi}$\textbf{.} The
sectors form a partition of $V$ into disjoint convex cones. Given sectors
$\mathcal{A}$ and $\mathcal{B}$, $\mathcal{A}$ is called a subsector of $B$ if
$A\subset B$.

A set $\mathcal{B}$ is called a face of $\mathcal{A}$ if its the intersection
of $C_{\phi}$ with a subspace $S_{\mathcal{T}}=\{y\in X:f(y-x)=0 $\textbf{,}
$\forall f\in\mathcal{T}\subset\overline{\mathcal{H}}\}$\textbf{.} In this
case we write $B\leq A$\textbf{.} A wall is a maximal proper face. These cells
will simply be referred to as conical cells based at $x$. \ Clearly if
$\mathcal{A}$ is a sector based at $x$, the set $\mathcal{A}-x+y$ is a sector
based at $y$\textbf{.}

A building $\Delta_{E}$ is called Euclidean if its \ apartments are Euclidean.
A geometric realization \ $\left\vert \Delta_{E}\right\vert $ of a Euclidean
Building is a building such that each apartment $\Sigma$ is a affine Coxeter complex.

\subsection{The Geometry of Euclidean Buildings\label{secaoCAT}}

Euclidean buildings have interesting geometric properties and can be thought
of as either an $n-$dimensional generalization of trees and some kind of
\ simplicial countenparts of symmetric space of non-positive curvature.

We start by defining a very special metric in the geometric realization
\ $\left\vert \Delta_{E}\right\vert $ of $\Delta_{E}$\textbf{.}Given two
points $x,y\in|\Delta_{E}|$ , the axiom $\mathbf{(B1)}$ of buildings says that
there exists a apartment $\Sigma$ containing $x$ and $y$. \ This apartment can
be endowed with a metric under which any chamber has diameter. Let $d_{\Sigma
}(x,y)$ denote the distance of $x$ and $y$ in this metric. The
\textbf{canonical metric} is then set to be:%
\[
d(x,y)=d_{\Sigma}(x,y)
\]

A curve $\gamma:I\rightarrow\left\vert \Delta_{E}\right\vert $ from the unit
interval $I$ to the metric space $\left\vert \Delta_{E}\right\vert $ is a
\textbf{geodesic} if there is a constant $\lambda$ $\geq0$ such that for any
$t\in I$ there is a neighborhood $J$ of $t$ such that for any $t_{1},t_{2}\in
J$ we have:%

\[
d(\gamma(t_{1}),\gamma(t_{2}))=\gamma|t_{1}-t_{2}|
\]

The \textbf{length} of a curve $\gamma$ is defined as
\[
\operatorname*{L}(\gamma)=\sup\left\{  \sum_{i=1}^{n}d(\gamma(t_{i}%
),\gamma(t_{i-1})):n\in\mathbb{N}\text{and } a=t_{0}<t_{1}<\cdots
<t_{n}=b\right\}  \text{.}%
\]
A r\textbf{ectifiable curve} is a curve with finite length.

In a metric space $(X,d)$, the induced \textbf{intrinsic metric}, $d_{i}(x,y)$
is the infimum of the lengths of all paths from $x$ to $y$. The length of such
a path is defined as explained for rectifiable curves. We set $d_{i}%
(x,y)=\infty$ if there is no path of finite length from $x$ to $y$. If
$d(x,y)=d_{i}(x,y)$ for all points $x$ and $y$ , we say $(X,d)$ is a length
space or a path metric space and the metric $d$ is intrinsic. The geometric
realization \ $\left\vert \Delta_{E}\right\vert $ with the canonical metric is
a length space.

\begin{proposition}
[\cite{Ga} pg. 194]Let a morphism $f:\Delta_{E}\rightarrow\Delta_{E}^{\prime}%
$\ \ of Euclidean buildings be given. Then $\left\vert f\right\vert
:\left\vert \Delta_{E}\right\vert \rightarrow\left\vert \Delta_{E}\right\vert
^{\prime}$ is an isometry between the corresponding geometric realizations.
\end{proposition}

The geometric realization $\left\vert \Delta_{E}\right\vert $ with the
canonical metric is a $CAT(0)$ space. This means that the curvature of such
spaces is smaller or equal to zero, that is, the triangles in such spaces are
thinner than in $\mathbb{R}^{2}$ in the following sense:

Given three points $x,y,z$ in a metric space $X$, a \textbf{comparison
triangle }$\bigtriangleup$ in $\mathbb{R}^{2}$ is a triangle with vertices
$\overline{x},\overline{y},\overline{z}$ such that%

\[
d(x,y)=\overline{d}(\overline{x},\overline{y}),\text{ }d(y,z)=\overline
{d}(\overline{y},\overline{z})\text{ and }d(x,z)=\overline{d}(\overline
{x},\overline{z})
\]

A metric space $X$ is called a\textbf{\ }$\mathbf{CAT(0)}$ space if for any
geodesical triangle $\bigtriangleup$ in $X$ and $\overline{\bigtriangleup}$ a
comparison triangle in $\mathbb{R}^{2}$\textbf{.} We have that $d(x,y)\leq
\overline{d}(\overline{x},\overline{y})$ for any $x,y$ $\in\bigtriangleup
$\textbf{.}

A complete $CAT(0)$ space is called a Hadamard space. The geometric
realization $\left\vert \Delta_{E}\right\vert $ of $\Delta_{E}$ is a Hadamard
space (see \cite{Ga} p. 197). \ From now on, a $CAT(0)$ space means a complete one.

The $CAT(0)$ inequality implies the following properties of the building
$\left\vert \Delta_{E}\right\vert $:

\begin{enumerate}
\item Any two points $x,y$ are joined by a unique geodesic segment, which
varies continuously in its endpoints $x,y$;

\item Metric balls are convex and contractible; in particular it is a simply
connected space and all of its higher homotopy groups are trivial;

\item A flat $\left\vert \Lambda\right\vert $ is a set that is isometric to
some $R^{n}$. Every maximal flat is an apartment and vice-versa. So we have a
geometric description of the apartments in the geometric realization
$\left\vert \Delta_{E}\right\vert $.

\item The apartments are convex and every geodesic is contained in an
apartment. \ A geodesic is called regular if it is contained in only one apartment.
\end{enumerate}

\subsubsection{The Building at the Infinite}

Given $X$ a locally compact $CAT(0)$ space,\ two geodesic rays $\sigma(t)$ and
$\gamma(t)$ are called asymptotic if there exists a constant $c$ such \ that
$d(\sigma(t),$ $\gamma(t))\leq c$ for all $t\geq0$. This define an equivalence
relation on the set of asymptotic rays we define \ the \textbf{ideal boundary}
$\partial_{\infty}X$ or simply $\partial_{\infty}X$ as the set of equivalence
classes of asymptotic rays. The union $X\cup\partial_{\infty}X$ \ is denoted
by $\overline{X}$. Given $x\in X$ and $\xi\in\partial_{\infty}X$, there exists
a geodesic ray $\gamma$ starting at $x$ and such that $\gamma(\infty)=\xi$.
This ray is denoted by $\gamma_{x,\xi}$

We introduce a topology at $\overline{X}=X\cup\partial_{\infty}X$ using as
base to this topology the open sets of $X$ and the following opens sets:
\[
U(x,\xi,R,\varepsilon):=\{z\in\overline{X}:z\notin B(x,R)\text{ and }d\left(
\sigma_{x,z}(R),\sigma_{x,\xi}(R)\right)  <\varepsilon\}
\]
with $x\in X$\textbf{,}$\xi\in\partial_{\infty}X$ and $R$, $\varepsilon$ real
positive numbers. This topology is knwon as the \textbf{Busemman topology} and
it does not depend on the choice of the base point, and it turns $\overline
{X}$ into a compact space with $X$ open and dense in $\overline{X}$, known as
\textbf{Busemman topology compactification}.

A fundamental fact of the topology we have just defined is\texttt{\ }that
given a isometry $\gamma$ of a complete $CAT(0)$ space, the natural extension
of $\gamma$ to $\overline{X}$ is a homeomorphism.

The compactification of the geometric realization of $\Delta_{E}$,
$\partial_{\infty}(\left\vert \Delta_{E}\right\vert )$\textbf{,} can be
endowed with a spherical building structure:

Given a conical cell $c \in\left\vert \Delta_{E}\right\vert $, it defines a
\textbf{simplex at infinity} (or \textbf{ideal simplex}) $c_{\infty}$ as the
set consisting of all equivalence classes $\gamma( \infty)$, where $\gamma$ is
a geodesic ray contained in the conical cell $c$. Given an ideal simplex
$\sigma$, we denote by $c_{\sigma, x} $ the unique conical cell based at $x$
such that $c_{\sigma, x} (\infty) =\sigma$

Given two ideal simplices $\sigma,\sigma^{\prime}\subset\partial_{\infty}X$,
we say that $\sigma^{\prime}$ is a \textbf{face} of $\sigma$ if the conical
cell $c_{\sigma^{\prime},x}$ is a face of $c_{\sigma,x}$ for some (and hence
all) $x\in X$. In that case we write $\sigma^{\prime}\leq\sigma$. This defines
an order relation that turns $\partial_{\infty}(\left\vert \Delta
_{E}\right\vert )$ a geometric realization of a spherical building. This is
stated in the next theorem and is the key to the study of control sets in this work.

\begin{theorem}
[\cite{Ga} p. 279]Consider the ideal simplex structure $\partial_{\infty
}(\left\vert \Delta_{E}\right\vert )$. The face relation is well defined and
the poset of ideal simplices of $\partial_{\infty}(\left\vert \Delta
_{E}\right\vert )$ is the geometric realization of the spherical building
associated with $\left\vert \Delta_{E}\right\vert $, i.e., $\partial_{\infty
}(\left\vert \Delta_{E}\right\vert )=\left\vert \Delta_{S}\right\vert $.
Moreover, the apartments of $\partial_{\infty}(\left\vert \Delta
_{E}\right\vert )$ are in bijection with the maximal system of apartments of
$\left\vert \Delta_{E}\right\vert $\textbf{.}
\end{theorem}

\section{Algebraic Groups over Local Fields}

\label{4}

Basic bibliographic reference for this section is \cite{Hu}. In the first
subsection we introduce only notation, taking for granted all definitions.

Let $k$ be a field with $\operatorname*{char}(k)=0$\textbf{.} An
\textbf{ultrametric norm} in a field $k$ is a function $\left\vert
\,\,\right\vert $ :$k\rightarrow\mathbb{Z}_{+}$\textbf{,} that satisfy

\begin{enumerate}
\item $\left\vert x\right\vert =0$ $\Leftrightarrow x=0$

\item $\left\vert xy\right\vert =\left\vert x\right\vert \left\vert
y\right\vert $ and $\left\vert 1\right\vert =1$

\item $\left\vert x+y\right\vert \leq\max\{\left\vert x\right\vert ,\left\vert
y\right\vert \}$
\end{enumerate}

The last inequality is called \textbf{ultrametric inequality} or
\textbf{strong triangular inequality}.

A field with a norm that satisfies the strong triangular inequality is said to
be an \textbf{ultrametric field} or\textbf{\ Archimedean field}. A field is
said to be complete if it is complete in the topology induced\ by the norm.

A \textbf{local field} is a complete, locally compact ultrametric field.

The $p$-adic field are the paradigmatic examples of local fields. In what
follows, we assume that $k$ denotes a local field.

A norm in a local field is induced by a valuation $v:k^{*}\mapsto\mathbb{Z}$,
in the sense that $\left|  x \right|  = e^{-v(x)}$. Let $A \subset k$ be the
subring of elements with $v(x)\geq0$ and let $\pi\in k$ an element such that
$v(\pi)=1$. The quotient field $A/\pi A$ is called the \textbf{residue field}
of $k$ with respect to the valuation $v$ (norm $\left|  \cdot\right|  $).

\subsection{The spherical building of a algebraic group}

Let $G$ be a linear algebraic group of dimension $n$ over a field $k$ and
$G^{0}$ the connected component of the identity element $Id\in G$. An
algebraic group $T$ over $k$ is a \textbf{torus} of rank $n$, if $T$ becomes
isomorphic to $(GL_{1})^{n}$ \ after the extension of the base field to the
algebraic closure $\mathbf{k}$. \ A torus is \textbf{split}\emph{\ } if it is
isomorphic over $k$ to $(GL_{1})^{n}$.

A subgroup $H\subset$ $G$ \ is called parabolic when $G/H$ is a complete
variety. A Borel subgroup $P$ is a minimal parabolic subgroup. We choose a
Borel subgroup $P$ which contains a maximal split torus $T$. Let
$\mathcal{N}_{sph}$ denote the normalizer of $T$ in $G$.

\begin{theorem}
[\cite{Hu},\cite{Borel}]The triple $(G,P,\mathcal{N}_{sph})$ is a
$\mathbf{BN}$-pair to which we can associate a spherical building $\Delta_{S}$.
\end{theorem}

The following remarkable results can be seen as consequences of the previous theorem

\begin{proposition}
\cite[p. 110]{Bro}(Bruhat Decomposition) Each double class $PgP$\textbf{,}
$g\in G$ can be written as $PwP$\textbf{,} $w\in W$\textbf{.} Also, we have
that mapping $w\rightarrow PwP$ is a bijection of the Weyl group onto the
double class $PgP $\textbf{,} with $g\in G$\textbf{.}
\end{proposition}

\begin{proposition}
\cite[p. 111]{Ga} If $gP^{I_{1}}g^{-1}=P^{I_{2}}$ then $I_{1}=I_{2}=I$ and
$g\in P^{I}$\textbf{.}
\end{proposition}

A decomposition of an algebraic group as a semidirect product $Z\rtimes
N$\textbf{,} with $Z$ reductive and $N$ unipotent is called a \textbf{Levi
decomposition} and we call $Z$ the \textbf{Levi factor} of $G$.

\begin{proposition}
(Levi decomposition) Let $N$ be the unipotent radical of $P$ and $Zbe$ the
centralizer of the torus $T$\textbf{.} Then
\[
P=Z\ltimes N
\]

\end{proposition}

\subsection{The Euclidean building of a algebraic group over a local field}

The spherical building described in the previous section depends only on the
algebraic structure of $G$. When the base field is a local one, it is possible
to assign to $G$ other $\mathbf{BN}$-pair, that yields an Euclidean building
for $G$.

\begin{theorem}
\cite[Bruhat-Tits Theorem]{TI1} Let $G$ be a \ simple and simply connected
group over a local field. Then we can associate to $G$ an Euclidean building
$\Delta_{E}$ through a $\mathbf{BN}$-pair, and the group $G$ acts as the group
of automorphism of $\Delta_{E}$.
\end{theorem}

We give a brief description of the building $\Delta_{E}$. Let us consider the
Weyl group $W_{0}$ of the spherical building $\Delta_{S}$. The Weyl group $W$
of $\Delta_{E}$ is an affine group that has $W_{0}$ as its linear part:
$W:=W_{0} \ltimes\mathbb{Z}^{n}$. It acts as a reflection group in
$\mathbb{R}^{n}$ and we denote by $\left|  \Sigma\right|  $ the (geometric
realization of the) simplicial structure determined by this action of $W$ on
$\mathbb{R}^{n}$. Then $\Sigma$ is a basic apartment and the structure of
$\Delta_{E}$ is fully determined when we state that the cardinality of
\[
\operatorname*{adjacent}(C,i):=\{C^{\prime} |C \sim_{i} C^{\prime} \},
\]
independs on the chamber $C$ and the adjacency type $i$ and is constant equal
to the cardinality of the residue field $A/ \pi A$ of $k$. In other words,
there is a prime $p$ such that any given cell has exactly $p$ $i$-adjacent
cells for any adjacency type $i$. Reciprocally, considering the description of
such a building $\Delta_{E}$, it is possible to prove that the group of
automorphism of $\Delta_{E}$ is isomorphic to $G$. Full details can be found
in \cite[chapter 5]{Bro}.

Using the Bruhat-Tits theorem we can give the following description of the
Euclidean building:

\begin{theorem}
[\cite{Bro}, p\'{a}g 163]The building $\Delta_{E}$ associated to a group $G$
with a Euclidean $\mathbf{BN}$-pair is isomorphic to the flag complex of the
incidence geometry consisting \ of the maximal bounded subgroups of $G$, where
two distinct such subgroups $P,Q$ are incident if and only if $P\cap Q$ is a
maximal subgroup of $P$.
\end{theorem}

The subgroups $P_{J^{\prime}}={\displaystyle\bigcup_{w\in W^{\prime}}} BwB$
\ of $G$ are called \textbf{paraholic subgroups}. They are open and compact
(in the non-Archimedian topology of $G$).

The spherical building at the infinity associated with this Euclidean building
coincides with the previous one, in the sense that $\partial_{\infty}\left(
\left\vert \Delta_{E}\right\vert \right)  =\left\vert \Delta_{S}\right\vert $.

\subsection{The action of $G$ in $\Delta_{E}$}

Assuming the short description of the euclidean building $\Delta_{E}$ given in
the previous section, we present now a more detailed description of the action
of $G$ on $\Delta_{E}$ and the structure it determines in the group $G$ itself.

We define the following subgroups of $G$:%

\begin{tabular}
[c]{lll}%
$\mathcal{N}_{aff}$ & $=$ & stabilizer of the apartment $\Sigma$\\
$B$ & $=$ & pointwise fixer of the chamber $C$\\
$T$ & $=$ & $\mathcal{N}_{aff}\cap B$\\
$P$ & $=$ & stabilizer of the chamber $C_{\infty}$\\
$\mathcal{N}_{sph}$ & $=$ & stabilizer of the apartment $\Sigma_{\infty}$\\
$Z$ & $=$ & $\mathcal{N}_{sph}\cap P$%
\end{tabular}

The subgroup $P$ is denominated \textbf{minimal parabolic }or\textbf{\ Borel}
\textbf{subgroup}, and the subgroup $B$ is the \textbf{parahoric minimal, }or
\textbf{Iwahori subgroup. }We remark that $(B,\mathcal{N}_{aff})$ is the Tits
system that originates the euclidean building and $(P,\mathcal{N}_{sph})$ is
the Tits system associated with the spherical building.

The stabilizer of the apartment $\Sigma\subset\Delta_{E}$ and the stabilizer
of the apartment $\Sigma_{\infty}\subset\Delta_{S}$ coincide, i.e.,
$\mathcal{N}_{aff}=\mathcal{N}_{sph}$ \cite[p. 259]{Ga}. So the $\mathbf{BN}$
decomposition of $G$ associated to both euclidean and spherical buildings, has
the same $\mathbf{N}$-factor and hence we write simply $\mathcal{N=}$
$\mathcal{N}_{aff}=\mathcal{N}_{sph}$.

We remark also that $W=\mathcal{N}/T$ is the affine Weyl group and that
$W_{0}=\mathcal{N}/Z$ is the spherical one.

If we denote by $\mathcal{N}_{trans}$ the subgroup of $\mathcal{N}$ that acts
on the apartment $\Sigma$ by translations, then the Levi factor $Z=\mathcal{N}%
\cap P$ of the minimal parabolic $P$ with respect to the apartment $A_{\infty
}$ is equal to the subgroup $\mathcal{N}_{trans}$ \cite[p. 200]{Ga}.
\label{levi toru}

A sector $\mathcal{A}$ with base point $s$ contained in the apartment $\Sigma$
determines a spherical chamber $\mathcal{A}^{\infty}\subset\left\vert
\Delta_{S}\right\vert =\partial_{\infty}(\left\vert \Delta_{E}\right\vert )$.
We define the semigroup $Z_{\mathcal{A}}$ as the subset of $Z=\mathcal{N}\cap
P$ that sends $s$ to a point in $\mathcal{A}$ and $N$ as the unipotent
subgroup of $G$ that stabilizes the chamber $\mathcal{A}^{\infty}$ of
$\partial_{\infty}(\left\vert \Delta_{E}\right\vert )$.

The group $G$ admits the following decompositions, that are analogous to the
Lie case.

\begin{theorem}
\cite{Bruhat}, \cite[p. 100]{Ron artigo} \strut

\begin{enumerate}
\item (Iwasawa Decomposition) $G=KZ_{\mathcal{A}}N$\textbf{,} and the double
classes $K\backslash G/N$ are in bijective correspondence with $Z$.

\item (Cartan Decomposition) $G=KZ_{\mathcal{A}}K$ and the double classes
$K\backslash G/K$ are in bijective correspondence with $Z_{\mathcal{A}}$.
\end{enumerate}
\end{theorem}

Being $G$ be a simple and simply connected group over a local field, the
quotient group of $G$ by any closed subgroup can be provided with a structure
of manifold over the local field, so in particular the flag manifolds are
manifolds over the local field (\cite[LG.4.10 and LG.4.11]{SErre} ).

Let $J$ \ be the set of generators of the spherical Weyl group $W_{0}$. Given
$I\subset J$ we define two special groups $W_{I}=\left\langle j_{i}:i\in
I\right\rangle $ and $W^{I}=\left\langle j_{k}:k\notin I\right\rangle $, and
associate to them the parabolic subgroups $P_{I}=PW_{I}P$ and $P^{I}=PW^{I}P$.
\ It is well known that $G=\left\langle P^{I},P_{I}\right\rangle $ and
$P^{I}\cap P_{I}=P$. Given the flag manifold $G/P^{I}$ and the projection
$G\overset{\pi}{\rightarrow}G/P^{I}$, we denote by $\dot{g}$ the lateral class
$gP^{I}\subset G/P^{I}$.

Let $N_{I}$ be the maximal unipotent subgroup of $P^{I}$ and let $\overline
{N}_{I}$ be the maximal unipotent subgroup of $P_{I}$\textbf{.} Let $\sigma\in
W_{0}$ the element with maximal lenght, i.e., the unique element of the Weyl
group $W_{0}$ that $\sigma$ maps a chamber to its opposite and define
$\overline{N}=\sigma N\sigma^{-1}$.

The manifold $\overline{N}_{I}\cdot\dot{e}$ can be identified with
$\overline{N}_{I}$ because $n_{2}n_{1}^{-1}\in P^{I}$ \ iff $n_{1}=n_{2}$ . We
will show the (well known in the real case) result which assures that
$\overline{N}_{I}\cdot\dot{e}$ is open and dense in the non-Archimedean
topology of $G/P^{I}$. Given $w\in W$ we define $\overline{N}^{w}%
:=\overline{N}\cap w\overline{N}w^{-1}$ and call each $\overline{N}^{w}$ a
\textbf{Bruhat cell}.

The function $gP\rightarrow gP^{I}$ will allow us to restric to the proof to
the case $G/P$ with $P$ a minimal parabolic subgroup.

\begin{theorem}
(Bruhat cellular decomposition) Then the mapping $\eta\rightarrow\eta\cdot
\dot{w}$ of $\overline{N}^{w}$ in $G/P$ is injective and $G/P$ is the disjoint
union of the cells \ $\overline{N}^{w}\dot{w}$\textbf{,} $w\in W$.
\end{theorem}

\begin{proof}
The Bruhat and Levi decompositions yield%
\[
G=\sigma G=\bigsqcup_{w\in W}\sigma PwP=\bigsqcup_{w\in W}\sigma NZ\sigma
^{-1}wP=\bigsqcup_{w\in W}\sigma N\sigma^{-1}wP=\bigsqcup_{w\in W}\overline
{N}wP \text{.}
\]

We will use the fact that the former union is disjoint to construct a
decomposition of $G/P$ into disjoint Bruhat cells.

Since $wPw^{-1}$ is the stabilizer of $wP$ in $G$, the stabilizer of $\dot
{w}=wP\in G/P$ in $\overline{N}$ is $\overline{N}\cap wPw^{-1}=\overline
{N}\cap wNw^{-1}$ and we can therefore identify $\overline{N}wP$ with
$\overline{N}/(\overline{N}\cap wNw^{-1})$.

Since $\overline{N}=(\overline{N}\cap w\overline{N}w^{-1})(\overline{N}\cap
wNw^{-1})$\textbf{,} it follows that $\overline{N}wP=\overline{N}^{w}wP$ can
be identified with $\overline{N}^{w}\cdot\dot{w}$. \ Thus the flag manifold
$G/P$ is the disjunct union of the Bruhat cells $\overline{N}^{w}\dot{w}$.
\end{proof}

The Bruhat cell $\overline{N}\cdot\dot{e}$ is open, since it has the same
dimension of$\ G/P$. Moreover, the Bruhat cells $\overline{N}^{w}\cdot\dot{w}$
have codimension at least one and there are finitely many such Bruhat cells,
so the cell $\overline{N}\cdot\dot{e}$ is also dense in $G/P$:

\begin{corollary}
The open cell $\overline{N}\cdot\dot{e}$ is open and dense in $G/P$.
\end{corollary}

Because of this last corollary, $\overline{N}_{I}\cdot\dot{e}$ is called the
\textbf{open Bruhat cell} of $G/P_{I}$.

Let $\Delta_{E}$ be the euclidean building associated with $G$ and $\left\vert
\Delta_{E}\right\vert $ be its geometric realization. As we mentioned, the
ideal boundary of $\left\vert \Delta_{E}\right\vert $ is the geometric
realization $\left\vert \Delta_{S}\right\vert $ of the spherical building.
Moreover, $\left\vert \Delta_{E}\right\vert $ is (open) and dense in the
Busemman compactification
\[
\overline{\left\vert \Delta_{E}\right\vert }=\left\vert \Delta_{E}\right\vert
\cup\left\vert \Delta_{S}\right\vert =\left\vert \Delta_{E}\right\vert
\cup\partial_{\infty}\left(  \left\vert \Delta_{E}\right\vert \right)
\text{.}%
\]
We define $\partial_{\infty}\left(  \Delta_{E}\right)  $ as the set of
barycenters of cells in $\left\vert \Delta_{S}\right\vert $. This is a formal
definition and clearly the set $\partial_{\infty}\left(  \Delta_{E}\right)  $
has a structure of spherical building. If we consider the barycenter of each
geometric cell in $\overline{\left\vert \Delta_{E}\right\vert }$, we find that
the barycenter of cells in $\left\vert \Delta_{S}\right\vert $ are the
acumullation points of barycenter of cells in $\left\vert \Delta
_{E}\right\vert $. Hence we can define the (non-trivial) \textbf{Busemman
topology} on the spherical building $\Delta_{S}$ as the topology defined in
$\partial_{\infty}\left(  \Delta_{E}\right)  $ as a subspace of $\overline
{\left\vert \Delta_{E}\right\vert }$.

We can associate to every simplex $S$ an parabolic group $P^{I}$ and a
parabolic type $I$\textbf{,} so \ that we can define a mapping:
\begin{align*}
\pi &  :\partial_{\infty}(\Delta_{E})\rightarrow\bigcup_{I\subset J}G/P^{I}\\
&  :S\mapsto P^{I}%
\end{align*}

The mapping $\pi$ can be decomposed in a family of mappings $\pi^{_{\Theta}}%
$\textbf{,} such that
\[
\pi^{_{\Theta}}:\partial_{\infty}^{I}(\Delta_{E})\rightarrow G/P^{I}%
\]
with $\partial_{\infty}^{I}(\Delta_{E})$ being the inverse image $\pi
^{-1}(G/P^{I})$. The simplices in $\partial_{\infty}^{I}(\Delta_{E})$ are
called the simplices of $I-$type in $\partial_{\infty}(\Delta_{E})$.

We note that $G$ acts transitively on $\partial_{\infty}^{\emptyset}\left(
\Delta_{E}\right)  $ (or equivalently $\Delta_{S}^{\emptyset}$) and hence may
be identified wit $G/P$ and so we can give $\partial_{\infty}^{\emptyset
}(\Delta_{E})$ the quotient topology, called the \textbf{ultrametric topology}.

\begin{proposition}
The topological spaces $\partial_{\infty}^{\emptyset}(\Delta_{E})$ with the
Busemman topology and $G/P$ with the ultrametric topology are
homeomorphic.\label{bruhateinfinito}
\end{proposition}

\begin{proof}
The group $G$ acts continuously and transitively in $\partial_{\infty
}^{\emptyset}(\Delta_{E})$ with $P$ stabilizing a chamber, so we have a
continuous bijection $\pi:$ $G/P\rightarrow$ $\partial_{\infty}^{\emptyset
}(\Delta_{E})$ However, since $\partial_{\infty}^{\emptyset}(\Delta_{E})$ is
compact, $\pi$\texttt{\ }must be a homeomorphism.
\end{proof}

The compactness of the flag manifold $G/P$ is a direct consequence of the
previous lemma (compare with \cite[p. 55]{Margulis}). At long last, we can
conclude that:

\begin{proposition}
Given a chamber $C$ \ in $\partial_{\infty}^{\emptyset}(\Delta_{E})$, $N$ be
the maximal unipotent group fixing $C$ and $\overline{N}$ the maximal
unipotent group fixing \ the opposite chamber $-C$, \ the orbit $\overline
{N}\cdot C$ is open and dense in $\partial_{\infty}^{\emptyset}(\Delta_{E})$.
\end{proposition}

\subsection{The action of $G$ in $\left\vert \Delta_{E}\right\vert $}

From here on, whenever needed, we will always assume (and it will be clear
from the context) as a basic assumptions that $G$ is a simple and simply
connected group over a local field and $\Delta_{E}$ its associated Euclidean building.

As we know (Section \ref{secaoCAT}), the geometric realization $\left\vert
\Delta_{E}\right\vert $ of the building $\Delta_{E}$ is a $CAT(0)$ space, in
which $G$ acts as an isometry group. We start stating some fundamentals
properties of isometries of $CAT(0)$ spaces.

The \textbf{displacement function} $d_{g}:X\rightarrow\mathbb{R}^{+}$ of an
isometry $g$ is defined by $d_{g}(x)=d_{g}\left(  g(x),x\right)  $. The
\textbf{translation length} of $g$ is the number $\left\vert g\right\vert
:=\inf\left\{  d_{g}(x):x\in X\right\}  $. The set of points $x$ such that
$d_{g}(x)=\left\vert g\right\vert $ is denoted by $\operatorname{Min}\left(
g\right)  $.

An isometry $g$ is called \textbf{semi-simple }if $\operatorname{Min}\left(
g\right)  $ is non\texttt{-}empty.

\begin{proposition}
\cite[chapter 2]{BGS} With the notation above established:

\begin{enumerate}
\item $\operatorname{Min}\left(  g\right)  $ is $g$ invariant;

\item If $h$ is a isometry of $X$, then $\left\vert g\right\vert =\left\vert
hgh^{-1}\right\vert $ and $\operatorname{Min}\left(  hgh^{-1}\right)  =
h\operatorname{Min}\left(  g\right)  $;

\item $\operatorname{Min}\left(  g\right)  $ is a closed convex set.
\end{enumerate}
\end{proposition}

The isometries can be classified as

\begin{enumerate}
\item \textbf{elliptic } if $g$ has a fixed point;

\item \textbf{hyperbolic }if $d_{g}$ attains a strictly positive minimum;

\item \textbf{parabolic }if $d_{g}$ does not attains its minimum (in other
words, if $\operatorname{Min}\left(  g\right)  =\emptyset$).
\end{enumerate}

\begin{proposition}
\cite[chapter 2]{BGS}

\begin{enumerate}
\item An isometry $g$ of $\left\vert \Delta_{E}\right\vert $ is hyperbolic if,
and only if, there is a geodesic $\gamma:\mathbb{R}\rightarrow\left\vert
\Delta_{E}\right\vert $ \ which is translated non-trivially by $g$, namely,
there is $a>0$ such that $g\cdot\gamma\left(  t\right)  =\gamma\left(
t+a\right)  $, for every $t\in\mathbb{R}$. Such a geodesic is called an
\textbf{axis} for $g$.

\item The axes of a hyperbolic isometry $g$ are parallel to each other and its
union is $\operatorname{Min}\left(  g\right)  $
\end{enumerate}
\end{proposition}

A hyperbolic isometry $h$ of an Euclidean building is said to be
\textbf{regular} if every axis of $h$ is contained in just one apartment.

\begin{lemma}
Let $\{\Sigma_{i}\}$ be a family of apartments\ that contains an axis $\gamma$
of $h$. Then $U:={\bigcup}$ $\Sigma_{i}$ and $I:={\bigcap}$ $\Sigma_{i}$ \ are
both invariants under $h$ in the sense that $hU=U$ and $hI=I$. In particular,
if $h$ is a regular hyperbolic isometry, i.e., every axis of $h$ is contained
in just one apartment, this apartment is invariant by $h$.
\end{lemma}

\begin{proof}
Let $\Sigma_{i}$ be an apartment containing the axis $\gamma$. Then
$h\Sigma_{i}$ is an apartment containing $\gamma$ since $h$ leaves $\gamma$ invariant.

The family $\{\Sigma_{i}\}$ is finite (in the case of a local field) and since
$h$ is a bijection, we see that $h$ acts as a permutations of the family
$\{\Sigma_{i}\}$. Consequently, we have
\[
hU={\bigcup}h\Sigma_{i}={\bigcup}\Sigma_{i}=U
\]
and%
\[
hI={\bigcap}h\Sigma_{i}={\bigcap}\Sigma_{i}=I
\]
whence they are both invariant under $h$.
\end{proof}

\begin{proposition}
Let $\left\{  \Sigma_{i}\right\}  _{i=1}^{k}$ be a collection of apartments in
$\left\vert \Delta_{E}\right\vert $, $\left\vert \Lambda\right\vert
:=\cap\Sigma_{i}$ be a flat in the euclidean building $\left\vert \Delta
_{E}\right\vert $, and $h$ a hyperbolic isometry of $\left\vert \Lambda
\right\vert $ that preserves the simplicial structure of $\left\vert
\Lambda\right\vert $. Then:

\begin{enumerate}
\item There exist a hyperbolic isometry $\tilde{h}$ of $\left\vert \Delta
_{E}\right\vert $ such that $\tilde{h}|_{\left\vert \Lambda\right\vert }=h$.

\item If $\left\vert \Lambda\right\vert $ is not a maximal flat (i.e, an
apartment) then any two such extensions $\tilde{h}_{1},\tilde{h}_{2}$ differ
by an element of $W_{\left\vert \Lambda\right\vert }$, i.e., $\tilde{h}%
_{1}(\tilde{h}_{2})^{-1}\in W_{\left\vert \Lambda\right\vert }$, where
$\Sigma$ is an apartment containing $\left\vert \Lambda\right\vert $ and
$W_{\left\vert \Lambda\right\vert }$ stand for the subgroup of $G$ that fixes
$\left\vert \Lambda\right\vert $ pointwise.

\item If $\left\vert \Lambda\right\vert $ is a maximal flat this extension is unique.
\end{enumerate}
\end{proposition}

\begin{proof}
Let $\Sigma$ be a maximal flat containing $\left\vert \Lambda\right\vert $.
The flat $\Sigma$ is isometric to $\mathbb{R}^{n}$, with $n$ the rank of the
building, and therefore $\left\vert \Lambda\right\vert $ is isometric to a
subspace of $\mathbb{R}^{n}$. We write $\Sigma=\left\vert \Lambda\right\vert
\oplus\left\vert \Lambda\right\vert ^{\perp}$, where $\left\vert
\Lambda\right\vert ^{\perp}$ is the orthogonal complement of $\left\vert
\Lambda\right\vert $ in $\Sigma$. Given a decomposition $x=x_{1}+x_{2}$ with
$x_{1}\in\left\vert \Lambda\right\vert $and $x_{2}\in\left\vert \Lambda
\right\vert ^{\perp}$, the map $\hat{h}$ $x=hx_{1}+x_{2}$ is an isometry of
$\Sigma$ that extends $h$ and preserves the simplicial structure of $\Sigma$.
Now to extend $\hat{h}\ $\ to $\left\vert \Delta_{E}\right\vert $, consider
$C$ a chamber in $\Sigma$ and $\hat{h}(C)$ its image. The second and third
axioms of buildings (Definition \ref{definicaoedificio}) assure the existence
of an isometry $\tilde{h}$ of $\left\vert \Delta_{E}\right\vert $ such that
$\tilde{h}(C)=\hat{h}(C)$ and $h(\Sigma$ $)=\Sigma$. This isometry satisfies
$\tilde{h}|_{\left\vert \Lambda\right\vert }=h$.

By definition, we have that $\tilde{h}_{1}(\tilde{h}_{2})^{-1}|_{\left\vert
\Lambda\right\vert }$ is the identity, and so, $\tilde{h}_{1}(\tilde{h}%
_{2})^{-1}\in W_{\left\vert \Lambda\right\vert }$. In particular, if
$\left\vert \Lambda\right\vert =\Sigma$ is a maximal flat we have that the
composition $\tilde{h}_{1}(\tilde{h}_{2})^{-1}$ fixes the apartment $\Sigma$
pointwise and the uniqueness is established.
\end{proof}

In a way analogous to the real Lie group case we have that
\[
H_{\Sigma}:=\left\{  h\in G:h\Sigma=\Sigma,h\text{ hyperbolic}\right\}
\]
is a maximal torus (page \pageref{levi toru}) and an isometry is regular if
and only if it is contained in a unique maximal torus. If we consider the
action of the Weyl group $W$ on an apartment we can give another
characterization of the regularity of $h$. Let $\gamma$ be an axis of $h$ and
$\Sigma$ an apartment containing $\gamma$. Given an hyperplane $P\subset
\Sigma$, we denote by $s_{P}$ the reflection in $P$ . Consider the set of
hyperplanes that determine reflections in $W$: $\mathcal{H}=\{P:s_{P}\in W\}$.
Considering the action of $W$ on $\Sigma$, an hyperbolic isometry $h$ is
\textbf{of type} $\Theta$ if $h$ has an axis invariant by a special subgroup
$W_{\Theta}$, with $W_{\Theta}$ a maximal special subgroup with this property.
The set $H_{\Sigma}^{\Theta}:=\left\{  h\in H_{\Sigma}:h\text{ is of type
}\Theta\right\}  $ determines a decomposition $H_{\Sigma}=\mathring{\cup
}_{\Theta\subset I}H_{\Sigma}^{\Theta}$ (disjoint union) with $H_{\Sigma
}^{\emptyset}$ having the dimension of $H_{\Sigma}$ and $\dim H_{\Sigma
}^{\Theta}<\dim H_{\Sigma}^{\Theta^{\prime}}$ whenever $\Theta^{\prime
}\varsubsetneqq\Theta$ (details can be found in \cite{Hu}). In a similar way,
given a flat $\left\vert \Lambda\right\vert $ in $\left\vert \Delta
_{E}\right\vert $ we define $H_{\left\vert \Lambda\right\vert }:=\left\{  h\in
G:h\left\vert \Lambda\right\vert =\left\vert \Lambda\right\vert ,h\text{
hyperbolic}\right\}  $

Considering any isometric identification of a flat $\left\vert \Lambda
\right\vert $ \ with the vector space $\mathbb{R}^{n}$, we fix $x_{0}%
\in\left\vert \Lambda\right\vert $ and associate to each hyperbolic isometry
$h$ a vector $v_{h}=h(x_{0})-x_{0}$, and denote $V_{H_{\left\vert
\Lambda\right\vert }}:=\{v_{h}:h\in H_{\left\vert \Lambda\right\vert }%
\}$\textbf{.} The action of $\mathbb{Z}$ in $V_{H_{\left\vert \Lambda
\right\vert }}$ defined by $nv_{h}=v_{h^{n}}$ turns $V_{H_{\left\vert
\Lambda\right\vert }}$ into a $\mathbb{Z-}$module. To each set of hyperbolic
isometries $K\subset H_{\left\vert \Lambda\right\vert }$, we associate a
vector space $V_{K}$ spanned by $\left\{  v_{h}:h\in K\right\}  $. Since the
action of $W$\ is irreducible, it follows that $\left\vert \Lambda\right\vert
$\ $=\left\langle v_{h}:h\in H_{\left\vert \Lambda\right\vert }\right\rangle $.

\bigskip

Up to the moment we studied the action of hyperbolic isometries on single
apartments or flats. Now we want to characterize the action of regular
hyperbolic isometries in the Euclidean building and in the spherical building.

We stress that, henceforth we consider only buildings associated to a simple
and simply connected group $G$ over a local field and hence, all isometries of
the associated buildings will be considered to be automorphisms of the
building (contained in $G$).

In all that follows, we use equal symbols to designate a simplex of the
Euclidean building $\left\vert \Delta_{E}\right\vert $\ or of its associated
spherical building $\partial_{\infty}(\left\vert \Delta_{E}\right\vert
)$,\ the sole distinction being that the latter are in boldface.

Let $h$ be a regular hyperbolic isometry, $\left\vert \Lambda\right\vert $ the
apartment invariant by $h$ and $\mathbf{\left\vert \Lambda\right\vert
}=\left\vert \Lambda\right\vert _{\infty}$ the apartment in the spherical
building $\partial_{\infty}(\left\vert \Delta_{E}\right\vert )$ that is the
boundary of the apartment $\left\vert \Lambda\right\vert $. The \textbf{main
attractor }of a hyperbolic isometry $h$ is the point $\xi\in\partial_{\infty
}(\left\vert \Delta\right\vert )$ such that $h^{n}(x)\rightarrow\xi$ for all
$x\in X$. If $h$ is regular, there is a unique open chamber at the infinite
containing $\xi$, denoted by $\mathbf{C}(\xi)$. \ The main attractor of
$h^{-1}$ is called the \textbf{principal repulsor }of $h$, denoted by $-\xi$
and $\mathbf{C}(-\xi)$ is the chamber opposite to $\mathbf{C}(\xi)$ in
$\mathbf{\left\vert \Lambda\right\vert }$\textbf{, } where $\left\vert
\Lambda\right\vert $ is the invariant apartment of $h$.

The main purpose of this section is to understand the behavior of the orbit of
an element $\eta\in\partial_{\infty}(\left\vert \Delta\right\vert )$ under the
action of a regular hyperbolic isometry $h$.

Our first observation is that when $\eta$ is in the apartment invariants by
$h$, i.e., $\eta\in\mathbf{\left\vert \Lambda\right\vert =}\partial_{\infty
}(\left\vert \Lambda\right\vert )$\textbf{,} the action of $h$ in $\eta$ is
trivial: $h\eta=\eta$. So we can restrict ourselves to the case $\eta
\notin\mathbf{\left\vert \Lambda\right\vert }$ . Next we consider an apartment
$\mathbf{\left\vert \Lambda\right\vert }^{\prime}$ (not necessarily unique)
containing $\eta$ and $-\xi$. In such a case, we have%
\[
h^{k}(\eta)\in\partial_{\infty}(h^{k}(\left\vert \Lambda\right\vert ^{\prime
}))\backslash\partial_{\infty}(h^{k}(\left\vert \Lambda\right\vert
))=\partial_{\infty}(h^{k}(\left\vert \Lambda\right\vert ^{\prime}%
))\backslash\partial_{\infty}(\left\vert \Lambda\right\vert )
\]
since $\left\vert \Lambda\right\vert $ is $h$-invariant.

To give a more accurate description of the orbit we will study the behavior of
some geodesic segments that converge to the points $\eta$ and $\xi$. In our
approach, we will focus on the Euclidian building, the corresponding results
for the spherical building being obtained by "projecting at infinity".

\begin{lemma}
[\cite{Bro}, p\'{a}g. 176]\label{lemaA}If two apartments $\left\vert
\Lambda\right\vert $ and $\left\vert \Lambda\right\vert ^{\prime}$ have a
common chamber $C$\texttt{\ }at infinity ($C\in\mathbf{\left\vert
\Lambda\right\vert }\cap\mathbf{\left\vert \Lambda\right\vert }^{\prime})$,
then $\left\vert \Lambda\right\vert \cap\left\vert \Lambda\right\vert
^{\prime}$ contains a sector $S$ such that $\partial_{\infty}(S)=\mathbf{C}$.
\end{lemma}

\begin{lemma}
\label{lemaB}Let $h$ be a regular hyperbolic isometry, and $\xi$ and $-\xi$
its main attractor and repulsor, respectively. Let $\left\vert \Lambda
\right\vert $ be the apartment invariant by $h$ and $S\subset\left\vert
\Lambda\right\vert $ a sector with $-\xi\in\partial_{\infty}(S)$. Then
$h^{k}S$ is a family of increasing nested sectors, in the sense that
$h^{k}(S)\varsubsetneq h^{k+1}(S)$. Moreover, given $x$ in the interior of $S$
and $\gamma$ a geodesic ray starting at $x$ and not entirely contained in $S$,
the length $L(\gamma|_{h^{k}(S)})$ of the segment of $\gamma$ contained in
$h^{k}(S)$ grows linearly.
\end{lemma}

\begin{proof}
Since $\left\vert \Lambda\right\vert $ is invariant by $h$ and $\left\vert
\Lambda\right\vert $ is isometric to $\mathbb{R}^{n}$, we translate the
problem to an euclidean setting in the which $S$ is a (simplicial) cone $h$
acts as a translation $hx=x+v_{h}$ where $-\lambda v_{h}\in S$ for every
$\lambda\geq\lambda_{0}$ for some $\lambda_{0}>0$ and $\gamma$ a staight line
passing throug points $y_{0}$ and $y_{1}$ with $y_{0}\in\operatorname{int}S$
and $y_{1}\notin S$. In this setting the problem becomes an elementary
euclidean question and it is easy to see that $h^{k}(S)\varsubsetneq
h^{k+1}(S)$ and $\operatorname{L}(\gamma|_{h^{k}(S)})=kd(x_{0},hx_{0})$ for
all $x_{0}\in\left\vert \Lambda\right\vert $\textbf{.}
\end{proof}

Given a point $x_{0}$ in the interior of $S$\textbf{,} an hyperbolic isometry
$h\in G$ and a point $\eta\in\partial_{\infty}\left(  \left\vert \Delta
_{E}\right\vert \right)  $ and an integer $k$ we define $\gamma^{k}(t)\ $as
the geodesic ray that starts at $x_{0}$ and such that $\gamma^{k}%
(\infty)=h^{k}(\eta)$.

Let $\sigma$ be the geodesic ray entirely contained in $\left\vert
\Lambda\right\vert $ with initial point $x_{0}$ and that is parallel to the
segment $\gamma^{0}|_{\left\vert \Lambda\right\vert }$. An initial segment of
$\sigma$ ($\left\{  \sigma\left(  t\right)  :0\leq t\leq\lambda\right\}  $ for
some $\lambda\geq0$) coincides with an initial segment of $\gamma^{0}$ and in
some point they fork as in a tree, with $\sigma$ being prolonged in
$\left\vert \Lambda\right\vert $ and $\gamma^{0}$ being prolonged in
$\left\vert \Lambda\right\vert ^{\prime}$. \ The same happens with $\sigma$
and $\gamma^{k}$, with the difference that $\gamma^{k}$ is prolonged in
$h^{k}(\left\vert \Lambda\right\vert ^{\prime})$. The lemma below describes
the growth of this intersection

\begin{lemma}
\label{lemaC}Let $\sigma$ and $\gamma^{k}$ be defined as above. Then, the
intersection of $\sigma$ and $\gamma^{k}$ is a geodesic segment whose length
grows linearly
\end{lemma}

\begin{proof}
By construction we know that $\sigma\cap\gamma^{0}=\sigma\cap S$. Now we prove
that $\sigma\cap\gamma^{1}=\sigma\cap h(S)$.

Consider the geodesic segment $h\gamma^{0}$. \ This segment passes through
$h(x_{0})$ and $h(\eta)$, and is contained in $h(\left\vert \Lambda\right\vert
^{\prime})$\textbf{.} Also, $h\gamma^{0}\cap\left\vert \Lambda\right\vert $ is
parallel to $\sigma$.

Now let $l$ be the geodesic contained in $h(\left\vert \Lambda\right\vert
^{\prime})$ that is parallel to $h\gamma^{0}$ and passes through $x_{0}$. The
segment $l\cap\left\vert \Lambda\right\vert $ is parallel to $\sigma$, i.e.,
$\sigma\cap\gamma^{1}$ $=\sigma\cap h(S)$. By the same argument we conclude
that $\sigma\cap\gamma^{k}=\sigma\cap h^{k}(S)$ and, by the previous lemma, we
find that the length of this intersection grows linearly, i.e.,
\[
\operatorname{L}(\sigma\cap\gamma^{k})=\operatorname{L}(\sigma\cap\gamma
^{0})+kd(x_{0},hx_{0}),\forall x_{0}\in\left\vert \Lambda\right\vert
\mathbf{.}%
\]

\end{proof}

Now we can state our main result describing the asymptotic action of
hyperbolic isometries.

\begin{proposition}
Let $h$ be a regular hyperbolic isometry, $\left\vert \Lambda\right\vert $ its
invariant flat and $\xi$ its main attractor. Let $G=KZN$ be the Iwasawa
decomposition defined by $\xi$ and $\left\vert \Lambda\right\vert $ ($N\xi
=\xi$ and $h\in Z$). Let $\overline{N}$ be the unipotent group that fixes
$-\xi$. Let $C(\eta)\in$ $\partial_{\infty}^{\emptyset}(\Delta_{E})$ be a
chamber and $C(\eta)=nw_{0}C(\xi)\in\partial_{\infty}^{\emptyset}(\Delta_{E}%
)$, with $w_{0}\in W$ and $n\in\overline{N}$ its expression relative to the
given Iwasawa decomposition. Then:
\[
\lim_{k\rightarrow\infty}h^{k}nw_{0}C(\xi){=}w_{0}C(\xi)
\]

\end{proposition}

\begin{proof}
Given $C(\eta)=nw_{0}C(\xi)$ \ and $-C(\xi)$, there is a flat $\left\vert
\Lambda\right\vert ^{\prime}\subset\left\vert \Delta_{E}\right\vert $ such
that $C(\eta),-C(\xi)\in\mathbf{\left\vert \Lambda\right\vert }^{\prime
}:=\partial_{\infty}\left(  \left\vert \Lambda\right\vert ^{\prime}\right)  $.
Let $\left\vert \Lambda\right\vert $ be the invariant flat of $h$ and
$\mathbf{\left\vert \Lambda\right\vert :=}$\textbf{$\partial_{\infty}$%
}$\left(  \left\vert \Lambda\right\vert \right)  $. Since $-\xi\in
\mathbf{\left\vert \Lambda\right\vert }\cap\mathbf{\left\vert \Lambda
\right\vert }^{\prime}$, Lemma \ref{lemaA} assures that the intersection of
$\left\vert \Lambda\right\vert $ and $\left\vert \Lambda\right\vert ^{\prime}$
contains a sector $S$ such that $\partial_{\infty}(S)=\mathbf{C}(-\xi)$.

Given $x_{0}\in\operatorname{int}S$ let $\sigma$ be the geodesic ray with
initial point at $x_{0}$ such that $\sigma\left(  \infty\right)  =nw_{0}\xi$
and for each integer $k$ define the $\gamma^{k}(t)$ as the geodesic ray with
$\gamma^{k}\left(  0\right)  =x_{0}$ and $\gamma^{k}\left(  \infty\right)
=h^{k}(\eta)$.

The rays $\gamma^{k}$ and $\sigma$ coincide in the intersection of the
apartments $\left\vert \Lambda\right\vert \cap h^{k}(\left\vert \Lambda
\right\vert ^{\prime})$ that increases linearly (by Lemma \ref{lemaC}). Then,
\cite[Proposition 8.19, pg. 268]{Bri} assures this characterizes the
convergence in the topology of Busemann, and we have that
\[
\lim_{k\rightarrow\infty}h^{k}nw_{0}\xi{=}w_{0}\xi
\]

\end{proof}

\begin{figure}[h]
\psfrag{e}{$\xi$} \psfrag{s}{$S$} \psfrag{n}{$\eta$} \psfrag{hn1}{$h\eta$}
\psfrag{F}{$\left\vert \Lambda \right\vert $}
\psfrag{f1}{$\left\vert \Lambda \right\vert '$}
\psfrag{Hf1}{$h\left\vert \Lambda \right\vert '$}
\includegraphics[height= 1.9464in,width= 4.8022in]{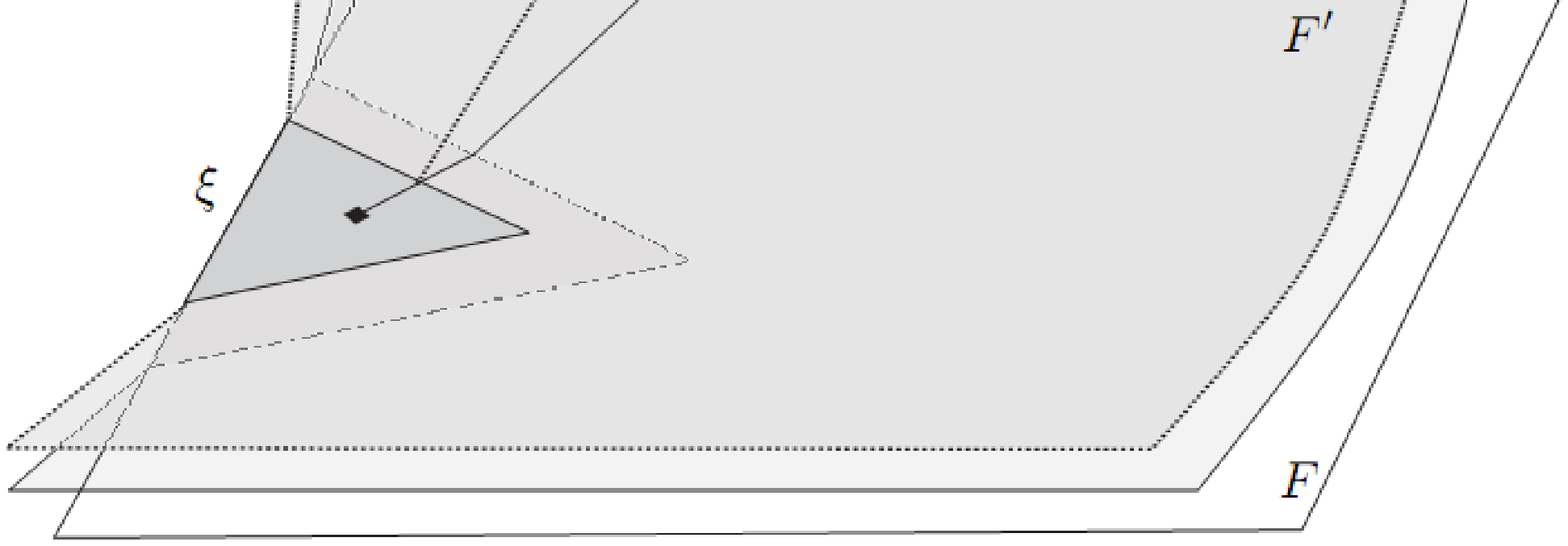}
\end{figure}

\section{Semigroups and Control Theory}

In this section, $G$ denotes a simple, simply connected algebraic group over a
local field. We remark that the group $G$ is also a Lie group over a local
field with the ultrametric topology $(G,\tau_{U})$\textbf{.} Our main object
of study are the semigroups $\mathcal{S}\subset G$ with non-empty interior in
the ultrametric topology.

Let $\xi\in\partial_{\infty}(\Delta_{E})$ be a regular point in the spherical
building and $C(\xi)$ the unique chamber that contains $\xi$. By Proposition
\ref{bruhateinfinito}, we can identify the chambers in $\partial_{\infty
}(\Delta_{E})$ with the points in the maximal flag $\mathbb{B}$ and so, by
abuse of notation, we write that $\mathbf{C}(\xi)\in$ $\mathbb{B}$, or that
$\xi\in\mathbb{B}$.

Given $\xi\in\mathbb{B}$, we denote by $B(\xi)$ the open Bruhat cell
containing $\xi$ .

\begin{proposition}
If \  a semigroup with non-empty interior $S\subset G$\ contains a hyperbolic
isometry of $\Delta_{E}$\textbf{,} then $\mathcal{S}$ contains a regular
hyperbolic isometry in its interior.
\end{proposition}

\begin{proof}
Since the simplicial structure of $\Delta_{E}$ determines a discrete metric,
the infimum of a translation function is attained, and so there are no
parabolic isometries in $G$. So we can assume $g\in\operatorname*{int}%
(\mathcal{S})$ as being a elliptic isometry. Let $x$ be a point fixed by $g$
and $K_{x}$ be the maximal compact subgroup fixing $x.($This subgroup is a
compact open subgroup of $G$! \cite{TI1} )  then $g\in\mathcal{S}%
_{x}=\operatorname*{int}(\mathcal{S})\cap K_{x},$ and so $\mathcal{S}_{x}$ is
a non-empty sub-semigroup of $K_{x}.$ Since $K_{x}$ is compact it follows that
$\mathcal{S}_{x}$ is a open subgroup of $K_{x}.$ Now, given a hyperbolic
isometry $h\in$ $\mathcal{S}$ since $g\in\mathcal{S}_{x}$ and $\mathcal{S}%
_{x}$ is a subgroup the isometry $ghg^{-1}$\ is hyperbolic in the interior of
$\mathcal{S}$\textbf{. }The existence of a regular hyperbolic isometry follows
then from the usual dimension arguments
\end{proof}

Unlike the real case, there are semigroups with non-empty interior that do not
contain any hyperbolic isometry.

\begin{proposition}
If  a semigroup with non-empty interior $S\subset G$\ contains nothing but
elliptic isometries, then $S$ is an open subgroup of $G$.
\end{proposition}

\begin{proof}
By the previous argument, $\mathcal{S}$ contains the inverse of any elliptic
isometry$g\in\operatorname*{int}(\mathcal{S})$\textbf{.} Therefore
$\operatorname*{int}(\mathcal{S})$\textbf{\ }is an open subgroup of $G$. If
$s\in\mathcal{S}$, and $g\in\operatorname*{int}(\mathcal{S})$\textbf{,} then
$gs\in\operatorname*{int}(\mathcal{S})$\texttt{.} But $g^{-1}\in
\operatorname*{int}(\mathcal{S})$ so that $s\in\operatorname*{int}%
(\mathcal{S})$.
\end{proof}

We assume henceforth that $\mathcal{S}\subset G$ is a semigroup of non-empty
interior that contains a hyperbolic isometry.

Let $h\in\operatorname{int}(\mathcal{S})$ be a regular hyperbolic isometry and
$\eta$ its main attractor.

\begin{lemma}
Let $\xi\in$ $\mathbb{B}$ and $B(\eta)$ the Bruhat cell based at $\eta,$
\ then exists  $g\in\mathcal{S}$ such that $g\xi\in B(\eta)$\textbf{.}
\end{lemma}

\begin{proof}
The proof is an immediate consequence from the fact that $B(\eta)$ is open and
dense on $\mathbb{B}$. Let $\xi\in\mathbb{B}$ \ and $A\subset
\operatorname{int}(\mathcal{S})$ be a open set. Then, $A\xi$ is open and
intercepts the Bruhat cell $B(\eta)$ that is dense in $\mathbb{B}$, i.e.,
$\exists g\in\operatorname{int}(\mathcal{S})$ such that $g\xi\in B(\eta)$.
\end{proof}

The flag manifold $\mathbb{B}$ is a compact space (\cite[p\'{a}g.
55]{Margulis}) and so the theorem \ref{existenciadecciemcompactos} yields the
existence of an invariant control set in $\mathbb{B}$ for the action of
$\mathcal{S}$. The following result assures the uniqueness of such a control
set in $\mathbb{B}$

\begin{theorem}
Let $\mathcal{S}\subset G$ be a semigroup with $\operatorname{int}%
(\mathcal{S})\neq\emptyset$ and that contains a hyperbolic isometry. Then the
invariant control set in $\mathbb{B}$ for the action of $\mathcal{S}$ is unique.
\end{theorem}

\begin{proof}
We will prove that there exists $\eta$ such that $\eta\in\operatorname{cl}%
(\mathcal{S}\xi)$ for all $\xi\in\mathbb{B}$, i.e., $\eta\in{\textstyle\bigcap
\limits_{\xi\in\mathbb{B}}}\operatorname{cl}(\mathcal{S}\xi)$. . \ Let $\eta$
be a main attractor for a regular hyperbolic isometry $h\in\operatorname{int}%
(\mathcal{S})$. By the previous lemma we can assume that $\xi\in B(\eta
)$\textbf{.} So we have that $\lim\limits_{k\rightarrow\infty}h^{k}(\xi)=\eta
$\textbf{,} for all $\xi\in B(\eta)$ and $\eta\in$ $\operatorname{cl}%
(\mathcal{S}\xi)$ for all $\xi$. Clearly this implies the uniqueness of the
invariant control set.
\end{proof}

We denote by $D$ the unique invariant control set in $\mathbb{B}$ and by
$D_{0}$ its set of transitivity.

We will now characterize the others control sets in $\mathbb{B}$.

We define%
\[
\Xi=\left\{  h\in\operatorname{int}(\mathcal{S}):h\text{ is a regular
hyperbolic isometry}\right\}
\]

Given a regular hyperbolic isometry $h$, there is a single apartment that
contains the axis of $h$ and a corresponding Iwasawa decomposition $G=KAN$.
The set of chambers in the spherical building fixed by $h$ is in bijection
with $W$, so we can associate\ to each of these chambers a $w$-type(which
depends on the choice of $h$) We will denote by $b(h,w)$ the point fixed by
$h$\ of $w$-type .

\begin{proposition}
Let $(D_{1})_{0}:=\{b(h,1):h\in\Xi\}$. Then $(D_{1})_{0}=D_{0}$%
.\label{caracterizacaoded0}
\end{proposition}

\begin{proof}
Since $(D_{1})_{0}$ is the set of main attractors we have that $(D_{1})_{0}$
is contained in $D_{0}$. To prove the other inclusion, we consider
$\eta=b(h,1)\in D_{0}$ with $h\in$ $\operatorname{int}(\mathcal{S})$ a regular
hyperbolic isometry. Given $\xi\in D_{0}$ exists $s_{1},s_{2}\in
\operatorname{int}(\mathcal{S})$ such that $s_{1}\eta=\xi$ and $s_{2}\xi=\eta
$, because $D_{0}$ is the set of transitivity of $D$. Now we consider
$h_{2}=s_{1}h^{n}s_{2}$. The following lemma proves that $s_{1}h^{n}s_{2}$is a
regular hyperbolic isometry for sufficiently large $n$.
\end{proof}

\begin{proposition}
If $\xi$\ is the main attractor of a regular hyperbolic isometry $h$ and
$s_{1},s_{2}$ are such that $s_{2}s_{1}$ belongs to the parabolic subgroup
that stabilizes $\xi$, then $s_{1}h^{n}s_{2}$\ is a regular hyperbolic
isometry for sufficiently large $n$ .
\end{proposition}

\begin{proof}
Let $g=s_{2}s_{1}$\textbf{,} then $s_{1}h^{n}s_{2}=s_{1}h^{n}gs_{1}^{-1}%
$\textbf{.} So it is enough to prove that $h^{n}g$ is hyperbolic when $g$ is
an isometry fixing the main attractor $\xi$ of $h$\textbf{.}

Let $\left\vert \Lambda\right\vert $ be the flat left invariant by $h$ and
$\left\vert \Lambda\right\vert ^{\prime}=g\left\vert \Lambda\right\vert $, $S$
be a sector in $\left\vert \Lambda\right\vert \cap\left\vert \Lambda
\right\vert ^{\prime}$ and $S_{2}=g^{-1}S\cap S$\textbf{,} so we have
$gS_{2}\subset S$\textbf{.}

Since $S_{2}$ and $gS_{2}$ are sectors in $\left\vert \Lambda\right\vert $
containing $\xi$\textbf{,} we have that, for $n$ large enough, $h^{n}%
gS_{2}\subset S_{2}$\textbf{,} i.e., $h^{n}g$ is a regular hyperbolic isometry.
\end{proof}

\begin{theorem}
For every $w\in W$ there exists a control set $D_{w}$ in $\mathbb{B}%
$\textbf{,}whose set of transitivity is%
\[
(D_{w})_{0}=\{b(h,w):h\in\Xi\}
\]

and these are all the control sets in $\mathbb{B}$.
\end{theorem}

\begin{proof}
The theorem will be proved in three steps.

In the first step we will prove that for any control set $D^{\prime}$ \ there
exists a fixed point of $w$-type $b(h,w)$ in $D^{\prime}$, for some regular
hyperbolic isometry $h$\textbf{.}

In the second step we will prove that, given two points with the same $w$-type
$\xi=b(h_{1},w)$ and $\eta=b(h_{2},w)$, with $w\in W$ and $h_{1},h_{2}\in\Xi
$\textbf{,} then $\xi\in\operatorname{cl}(\mathcal{S}\eta)$ and $\eta
\in\operatorname{cl}(\mathcal{S}\xi)$. Therefore, by Proposition
\ref{pontosfixosestaoem cc}, $(D_{w})_{0}$ is contained in the transitivity
set of a control set $D_{w}$.

At last, we will prove that $(D_{w})_{0}:=\{b(h,w):h\in\Xi\}$ is the set of
transitivity of $D_{w}$.

(1) Given a control set $D^{\prime}$ and $\eta\in D^{\prime}$, let $P$ be the
isotropy group of $\eta$. If $\eta\in D_{0}^{\prime}$ then $P\cap
\operatorname{int}(\mathcal{S})\neq\emptyset$. The subgroup $P$ admits an
Iwasawa decompositions as $P=MAN^{+}$\textbf{.} The subset
\[
\sigma=\{m\in M:\exists hn\in AN^{+}\text{with }mhn\in\operatorname{int}%
(\mathcal{S})\}
\]
has non-empty interior in $M$\textbf{,} since $M$ normalizes $AN$. The fact
that $M$ is compact implies that $\sigma$ is a subgroup of $M$. So
$\operatorname{int}(\mathcal{S})\cap AN^{+}\neq\emptyset$\textbf{,} and hence
there exists a regular hyperbolic isometry $g\in\operatorname{int}%
(\mathcal{S})$ such that $g\eta=\eta$. In particular, $\eta$ is is a fixed
point of $w$-type for $g$.

(2) Given $\xi=b(h_{1},w)$ and $\eta=b(h_{2},w)$, if $b(h_{1},1)$
$=b(h_{2},1)$ then $b(h_{1},w)$ is in the same Bruhat cell that as
$b(h_{2},w)$, so $\lim_{n\rightarrow\infty}(h_{2})^{n}b(h_{1},w)=b(h_{2},w)$.
If $b(h_{1},1)\neq b(h_{2},1)$\textbf{,} the Proposition
\ref{caracterizacaoded0} assures that exists $s_{1}\in\mathcal{S}$\textbf{,}
such that $s_{1}b(h_{1},1)=b(h_{2},1)$ since $(D_{1})_{0}$ is a transitivity
set, and thereby we reduce to the previous case.

(3) Now we will prove that for any control set $D^{\prime}$ its set of
transitivity is $(D_{w})_{0}$ for some $w$. Let $\eta=b(g,w)\in D^{\prime}$ be
the attractor of $w$-type in $D^{\prime}$ whose existence was proved in (1).

Given $\xi\in D_{0}^{\prime}$, we have $s_{1},s_{2}\in\operatorname{int}%
(\mathcal{S})$ such that $s_{1}\eta=\xi$ and $s_{2}\xi=\eta$. Then
$h_{2}=s_{1}g^{n}s_{2}$ is a regular hyperbolic isometry and $\xi=b($
$s_{1}h^{n}s_{2},w_{1})$. Thus $D_{w}$ is the set of transitivity of
$D^{\prime}$.
\end{proof}

\subsection{Subgroup $W(S)$}

In the previous section, we proved that any control set for $\mathcal{S}$ in
$\mathbb{B}$ is of the form $D_{w}$ for some $w\in W$ .

We now turn to the problem of determining when two such control sets
$D_{w},D_{w^{\prime}}$\ coincide. To start with, define a subset $W(S)$\ of
$W$\ by the rule%
\[
W(\mathcal{S})=\left\{  w\in W:D_{w}=D_{1}\right\}
\]

By its definition\texttt{,}$W(S)$\texttt{\ }depends on the choice of a chamber
$\mathbf{C}^{+}$ in the apartment $\mathbf{\left\vert \Lambda\right\vert
\subset}$\textbf{$\partial_{\infty}$}$_{\infty}(X)$ which we fix once and for
all. (Notice that if we had chosen the chamber $\mathbf{C}_{1}^{+}%
=g\mathbf{C}^{+}$\ instead of $\mathbf{C}^{+}$, then\texttt{ }%
\[
w\in W(\mathcal{S},\mathbf{C}^{+})\Leftrightarrow gwg^{-1}\in W(\mathcal{S}%
,g\mathbf{C}^{+})\mathbf{.}%
\]
so that conjugation by $g$\ defines an isomorphism between\texttt{
}$W(S,\mathbf{C}^{+})$\texttt{\ }and\texttt{ }$W(S,g\mathbf{C}^{+})$\texttt{,
}which does not depend on the choice of $g$\ taking\texttt{ }$\mathbf{C}^{+}%
$\texttt{\ }to\texttt{ }$g\mathbf{C}^{+}$\texttt{\ }since any two such
elements differ by an element that fixes pointwise the chamber\texttt{
}$\mathbf{C}^{+}$\textbf{.}

From now on, let\texttt{ }$b_{0}$\texttt{\ }be the image\texttt{ }in\texttt{
}$\in\mathbb{B}$\texttt{ }of the chamber $\mathbf{C}^{+}$ \texttt{(}used to
define\texttt{ }$W(S)=W(S,C^{+})$\texttt{)}.

\begin{lemma}
\label{LemaW}Given $b_{0}\in(D_{1})_{0}$ Then are equivalent:

\begin{enumerate}
\item $w\in W(\mathcal{S});$

\item $\tilde{w}b_{0}\in(D_{1})_{0}$\textbf{,} with $\tilde{w}$ being a
representative of $w$ in $M^{\ast}$, and $W=M^{\ast}/M$\textbf{.}
\end{enumerate}
\end{lemma}

\begin{proof}
Let $b_{0}\in(D_{1})_{0}$, $w\in W(\mathcal{S})$ and $\tilde{w}$ be a
representative of $w$ in $M^{\ast}$.By the characterization of $(D_{w}%
)$\textbf{,} $\tilde{w}b_{0}\in(D_{w})_{0}$\textbf{,} i.e., $\tilde{w}b_{0}%
\in(D_{w})_{0}=(D_{1})_{0}$. As for the converse, by definition any $w\in
W(S)$\ is such that $D_{w}=D_{1}$\textbf{,} hence $(D_{w})_{0}=(D_{1})_{0}%
$\textbf{. }However $(D_{w})_{0}=\{b(h,w):h\in\Xi\}$ and therefore $wb_{0}%
\in(D_{w})_{0}$\textbf{.}
\end{proof}

\begin{proposition}
$W(\mathcal{S})$ is a subgroup of $W$\textbf{.}
\end{proposition}

\begin{proof}
Let $b\in(D_{1})_{0}$, $w_{1}$ and $w_{2}\in W(\mathcal{S})$ and $\tilde
{w}_{1}$ and $\tilde{w}_{2}$ they representatives in $M^{\ast}$. By the Lemma
\ref{LemaW} $\tilde{w}_{1}b\in(D_{1})_{0}$. Another application of \ref{LemaW}
\ for $\tilde{w}_{1}b$ allows us to conclude that $(\tilde{w}_{1}\tilde{w}%
_{2}(\tilde{w}_{1})^{-1})\tilde{w}_{1}b\in D_{0}$ and therefore $\tilde{w}%
_{1}\tilde{w}_{2}b\in D_{0}$. Hence $W(\mathcal{S})$ is a semigroup of $W$,
but since $W$ is finite $W(\mathcal{S})$ is a subgroup.
\end{proof}

$W(\mathcal{S})$ is not only a subgroup of $W$\textbf{,} but a Weyl subgroup
of $W$, i.e., $W(\mathcal{S})=W_{\Theta}$ for some subset $\Theta$ of simple roots.

\begin{theorem}
$W(\mathcal{S})=W_{\Theta}$\textbf{,} for some $\Theta\subset\Pi$
\end{theorem}

\begin{proof}
Let $\mathbb{H}=\{h\in\operatorname{int}(\mathcal{S}):W(\mathcal{S})h=h\}$.

Suppose that $\mathbb{H}\neq\emptyset$ and let $\widehat{h}\in\mathbb{H}$ be
an isometry of maximum regularity in $\mathbb{H}$\textbf{,} \ and $\Theta$ be
its type. We will show that $W(\mathcal{S})=W_{\Theta}$ proving that
$W(\mathcal{S})$ acts transitively in $\mathcal{C}_{h}=\{$the chambers that
have $\widehat{h}$ as a wall$\}$.

Since $\widehat{h}\in\operatorname{int}(\mathcal{S})$ there is a regular
hyperbolic isometry $h_{\mathbf{C}}\in\operatorname{int}(\mathcal{S})$ whose
main attractor is the chamber $\mathbf{C}\in\mathcal{C}_{h}$.

So every chamber $\mathbf{C}\in\mathcal{C}_{h}$ belongs to the set of
transitivity of the invariant control set in $\mathbb{B}$ and so $\mathcal{S}
$ acts transitively in $\mathcal{C}_{h}$.

Suppose now that $H=\emptyset$; we will prove that in this case
$W=W(\mathcal{S})$. Let $\mathbf{C}$ be a chamber in the set of transitivity
$D_{0}$ of the invariant control set in $\mathbb{B}$ and let $h$ be a regular
hyperbolic isometry that has $\mathbf{C}$ as its the main attractor.

We denote $h_{w}=wh^{n}$ and remark that all $h_{w}$ are regular hyperbolic
isometries for $n$ sufficiently large. Let $K$ \ be the cone and
$K_{\mathbb{N}}$ be the lattice defined as follow:%
\begin{align*}
K  &  =\{\sum_{w\in W}a_{w}v_{h_{w}}:a_{w}\in\mathbb{R}^{+}\mathbb{\}}\\
K_{\mathbb{N}}  &  =\{\sum_{w\in W}a_{w}v_{v}h_{w}:a_{w}\in\mathbb{N\}}%
\end{align*}

Notice that $K$ is a vector space, because the absence of fixed points implies
that $v_{h^{1}}+\ldots+v_{h^{w}}=0$. And as the action of $W$ is irreducible
and $W(\mathcal{S})\subset W$ acts without fixed points we have that $K$t is
the whole apartment $\left\vert \Lambda\right\vert $.

Since $K_{\mathbb{N}}\subset\operatorname{int}(\mathcal{S})$ is a lattice in
$K$ we know that for any chamber $\mathbf{C}\in\partial_{\infty}(\left\vert
\Lambda\right\vert )$ there exists a hyperbolic isometry $h^{\prime}\in
K_{\mathbb{N}}\subset\operatorname{int}(\mathcal{S})$\textbf{,} with main
attractor $\mathbf{C}$\texttt{. }Therefore\texttt{ }$W(\mathcal{S})$ acts
transitively in the chambers of $\partial_{\infty}(\left\vert \Lambda
\right\vert )$. But, as $W$ acts simply transitive in the chambers of
$\partial_{\infty}(\left\vert \Lambda\right\vert )$\textbf{,} we conclude that
$W(\mathcal{S})=W$.
\end{proof}

The map $w\rightarrow D_{w}$ defined in the theorem of characterization of
($D_{w})_{0}$ it is not necessarily bijective. The following theorem says that
we can parametrize the control sets by the lateral classes $W(\mathcal{S}%
)\backslash W$.

\begin{theorem}
$D_{w_{1}}=D_{w_{2}}$ if and only if $w_{1}w_{2}^{-1}\in W(\mathcal{S})$.
Hence the control sets of $\mathcal{S}$ in $\mathbb{B}$ are in bijection with
$W(\mathcal{S})\backslash W$.
\end{theorem}

\begin{proof}
It suffices to show that $(D_{w_{1}})_{0}=(D_{w_{2}})_{0}$ if and only if
$w_{1}w_{2}^{-1}\in W(\mathcal{S})$.

($\Rightarrow)$ Let $b_{0}\in D_{0}=(D_{1})_{0}$. By the characterization of
the control sets, $\tilde{w}_{1}b_{0}\in(D_{1})_{0}$, and by hypothesis we
have $(D_{w_{1}})_{0}=(D_{w_{2}})_{0}$\textbf{. }Since $\tilde{w}_{1}%
b\in(D_{w_{2}})$ we have by the definition of $(D_{w_{2}})_{0}$ and by
conjugation, that $(\tilde{w}_{1}\tilde{w}_{2}(\tilde{w}_{1})^{-1})^{-1}%
\tilde{w}_{1}b\in D_{0}$, and so $(\tilde{w}_{2})^{-1}\tilde{w}_{1}\in
W(\mathcal{S})$.

($\Leftarrow$) Let $b^{\prime}=\tilde{w}_{1}b_{0}\in(D_{1})_{0}$ . Then
$(\tilde{w}_{2})^{-1}b^{\prime}\in D_{0}$, since $w_{2}^{-1}w_{1}\in
W(\mathcal{S})$. And $(\tilde{w}_{1}\tilde{w}_{2}(\tilde{w}_{1})^{-1}%
)^{-1}\tilde{w}_{1}b_{0}$ $\ \in D_{0}$ and consequently $b^{\prime}%
\in(D_{w_{2}})_{0}$\textbf{,} by the definition of $(D_{w_{2}})_{0}$.
\end{proof}

\end{document}